\documentclass[12pt]{amsart}
\usepackage[colorlinks=true,citecolor=black,linkcolor=black,urlcolor=blue]{hyperref}
\usepackage{amsmath,pifont,enumitem,mathtools}
\usepackage[english,  activeacute]{babel}
\usepackage[utf8]{inputenc}
\usepackage[dvipsnames]{xcolor}
\usepackage{amssymb}
\usepackage{amsthm}
\usepackage{comment}
\usepackage{graphics,graphicx}
\usepackage{array}
\usepackage{bm}
\usepackage{a4wide}
\usepackage{color, url}
\usepackage{float}
\usepackage{xcolor}

\newcommand{\seqnum}[1]{\href{http://oeis.org/#1}{\underline{#1}}}

\usepackage[colorinlistoftodos]{todonotes}

\theoremstyle{plain}
\newtheorem{theorem}{Theorem}[section]

\newtheorem{coro}[theorem]{Corollary}

\theoremstyle{definition}

\newtheorem{example}[theorem]{Example}

\newtheorem{remark}[theorem]{Remark}
\newcommand{\Cat}{{\mathcal{C}}}

\def\val{{\textsf{val}}}
\def\wruns{{\textsf{wruns}}}
\def\wdruns{{\overline{\textsf{wruns}}}}

\def\lval{{\textsf{$\ell$-val}}}
\def\lpea{{\textsf{$\ell$-peak}}}
\def\peak{{\textsf{peak}}}
\def\runs{{\textsf{runs}}}
\def\druns{{\overline{\textsf{runs}}}}
\def\vsym{{\textsf{symv}}}
\def\psym{{\textsf{symp}}}
\def\wfwn{{f}}

\newcommand{\tx}[1]{\texttt{#1}}

\newcommand{\f}{{\textsf{Flat}}}

\title{Flattened Catalan Words}

\date{\today}

\begin{document}
\author[J.-L Baril]{Jean-Luc Baril}
\address[J.-L Baril]{LIB, Universit\'e de Bourgogne Franche-Comt\'e,   B.P. 47 870, 21078, Dijon Cedex, France}
\email{\textcolor{blue}{\href{mailto:barjl@u-bourgogne.fr}{barjl@u-bourgogne.fr}}}

\author[P.~E.~Harris]{Pamela E. Harris}
\address[P.~E.~Harris]{Department of Mathematical Sciences, University of Wisconsin-Milwaukee, Milwaukee, WI 53211 United States}
\email{\textcolor{blue}{\href{mailto:peharris@uwm.edu}{peharris@uwm.edu}}}

\author[J. L. Ram\'{\i}rez]{Jos\'e L. Ram\'{\i}rez}
\address[J. L. Ram\'{\i}rez]{Departamento de Matem\'aticas,  Universidad Nacional de Colombia,  Bogot\'a, Colombia}
\email{\textcolor{blue}{\href{mailto:jlramirezr@unal.edu.co}{jlramirezr@unal.edu.co}}}

\begin{abstract}
In this work, we define flattened Catalan words as Catalan words whose runs of weak ascents have leading terms that appear in weakly increasing order. We provide generating functions, formulas, and asymptotic  expressions for the number of flattened Catalan words based on the number of runs of ascents (descents), runs of weak ascents (descents), $\ell$-valleys, valleys, symmetric valleys, $\ell$-peaks, peaks, and symmetric peaks. 
\end{abstract}

\maketitle

\section{Introduction}
A word $w=w_1w_2\cdots w_n$ over the set of nonnegative integers is called a \emph{Catalan word} if $w_1=\tx{0}$ and $\tx{0}\leq w_i\leq w_{i-1}+\tx{1}$ for $i=2, \dots, n$.  Throughout this paper, $|w|$ denotes the length of~$w$ and $\epsilon$  denotes the \emph{empty word}, which is the unique word of length zero. For $n\geq 0$, let $\Cat_n$ denote the set of Catalan words of length $n$. 
We set  $\Cat\coloneq\bigcup _{n\geq0}\Cat_n$ and $\Cat^+\coloneq\bigcup _{n\geq1}\Cat_n$ be the set of nonempty Catalan words. 
For example,
\begin{align*}
\Cat_4=\left\{\begin{matrix} \texttt{0000}, \, \texttt{0001}, \, \texttt{0010}, \, \texttt{0011}, \, \texttt{0012}, \, \texttt{0100},  \, \texttt{0101},  \\\texttt{0110},   \, \texttt{0111}, \, \texttt{0112}, \,  \texttt{0120},  \,\texttt{0121}, \, \texttt{0122},  \,\texttt{0123}\end{matrix}\right\}.
\end{align*}
Note that $|\Cat_n|=c_n=\frac{1}{n+1}\binom{2n}{n}$ is the $n$th Catalan number.  
 The exploration of Catalan words has begun with the comprehensive generation of Gray codes tailored for growth-constricted words \cite{ManVaj}. Baril et al.  \cite{BGR, Baril2, Baril} have delved into analyzing the distribution of descents and the ultimate symbol in Catalan words avoiding one or two classical patterns of length at most three. Similar findings \cite{Baril3, Baril4, AlejaRam} emerge in studies of restricted Catalan words avoiding consecutive patterns of length three or pairs of relations. Callan et al. \cite{CallManRam} initiate the enumeration of statistics, including area and perimeter, on the polyominoes associated with Catalan words. Furthermore, assorted combinatorial statistics regarding polyominoes associated with both Catalan and Motzkin terminologies have been scrutinized \cite{Baril5, ManRamF, ManRamM, Toc}. Next Shattuck \cite{Shattuck} initiated an examination into the frequency of distinct subword occurrences, spanning no more than three characters, nestled within Catalan words,  like descents, ascents, and levels.
In a recent paper \cite{Paper1}, Baril et al.~provide generating functions, formulas, and asymptotic  expressions for the number of Catalan words based on the number of runs of ascents (descents), runs of weak ascents (descents), $\ell$-valleys, valleys, symmetric valleys, $\ell$-peaks, peaks, and symmetric peaks. 

Given a permutation of $[n]=\{1,2,\ldots,n\}$ in one-line notation $\pi=\pi_1\pi_2\cdots\pi_n$, the \emph{runs} of $\pi$ are the maximal contiguous increasing subwords of~$\pi$. If the sequence of leading terms of the runs of $\pi$ appears in increasing order, 
then $\pi$ is called \textit{flattened partition} of length $n$. 
Nabawanda et al.\ give recursive formula for the number of flattened partitions of length $n$ with $k$ runs~\cite[Theorem 1]{ONFRAB}.
Callan gives the number of flattened partitions of length $n$ avoiding a single 3-letter pattern~\cite{Callan}.  Elder et al.\ extended the work Nabawanda et al.\ to establish recursive formulas for the number of flattened parking functions built from permutations of $[n]$, with $r$ additional ones inserted that have $k$ runs~\cite[Theorems 29, 30 and 35]{flat_pf}.
A further generalization includes the work of Buck et al.\ \cite{FlatStirling} who establish that flattened Stirling permutations are enumerated by the Dowling numbers, which corresponds to the OEIS entry \cite[\seqnum{A007405}]{OEIS}. 

In this work, we define  \emph{flattened Catalan words}, which are Catalan words whose maximal contiguous nondecreasing subwords  have leading terms in weakly increasing order. For example, the Catalan word $\tx{0012301222345523343}\in \Cat_{19}$ is a flattened Catalan word with four maximal contiguous nondecreasing subwords $\tx{00123}$, $\tx{012223455}$, $\tx{2334}$, and $\tx{3}$, whose leading terms satisfy $\tx{0}\leq \tx{0}\leq \tx{2}\leq \tx{3}$. Conversely, $\tx{012321}\in\Cat_6$ is not a flattened Catalan word as it 
has maximal contiguous nondecreasing subwords $\tx{0123}$, $\tx{2}$, and  $\tx{1}$, and the leading terms $\tx{0}$, $\tx{2}$, and $\tx{1}$ are not in weakly increasing order. We denote the sets of nonempty flattened Catalan words and flattened Catalan words of length $n$ as $\f(\Cat^+)$ and $\f(\Cat_n)$, respectively.

Let $w=w_1w_2\cdots w_n\in \f(\Cat_n)$.
As usual, we say that $w$ has an \emph{ascent} (\emph{descent}) at position $\ell$ if $w_\ell < w_{\ell+1}$ ($w_\ell > w_{\ell+1}$), where $\ell \in [n-1]$. Similarly, we define \emph{weak ascent} (resp. \emph{weak descent})
at position $\ell$ if $w_\ell \leq w_{\ell+1}$ ($w_\ell \geq w_{\ell+1}$), where $\ell \in [n-1]$.
A \emph{run} (resp. \emph{weak run}) of ascents (resp. \emph{weak ascents}) in a word $w$ is a maximal subword of consecutive ascents (resp. weak ascents). 
The number of runs  in $w$ is denoted by $\runs(w)$, and the number of weak runs  in $w$ is denoted by $\wruns(w)$. 
The runs of descents and weak descents are defined similarly, and the  statistics will be  denoted $\druns(w)$ and $\wdruns(w)$, respectively. 
An $\ell$-\emph{valley}  in a flattened Catalan word $w$ is a subword of the form $ab^\ell(b+1)$, where $a>b$ and  $\ell$ is a positive integer and $b^\ell$ denotes $\ell$ consecutive copies of the letter $b$.  
If $\ell=1$, we say that it is a  \emph{short valley}.  
The number of $\ell$-valleys of $w$ is denoted by $\lval(w)$ and the number of  all $\ell$-valleys for $\ell\geq 1$ of $w$ is denoted by $\val(w)$.  
A \emph{symmetric valley} is a valley of the form $a(a-1)^\ell a$ with $\ell\geq 1$. The number of symmetric valleys of $w$ is denoted by $\vsym(w)$. Analogously, we define the peak statistic. Namely, an $\ell$-\emph{peak} in $w$ is a subword of the form $a(a+1)^\ell b$, where $a\geq b$ and  $\ell$ is a positive integer.
The number of $\ell$-peaks of $w$ is denoted by $\lpea(w)$ and
the sum of all $\ell$-peaks for $\ell\geq 1$ of $w$ is denoted by $\peak(w)$. 
If $\ell=1$, we say that it is a  \emph{short peak}; and if $a=b$,  it is called a \emph{symmetric peak}. The number of symmetric peaks of $w$ is denoted by $\psym(w)$.  

Our contributions include generating functions and combinatorial expressions for the number of flattened Catalan words based on the number of 
runs of ascents (descents), runs of weak ascents (descent), $\ell$-valleys, valleys, symmetric valleys, $\ell$-peaks, peaks, and symmetric peaks.  We also establish one-to-one correspondences between:
\begin{itemize}
    \item flattened Catalan words of length $n$ with $k$ runs of ascents and $k$-part order-consecutive partitions of $n$, which have been studied in \cite{WM}, see Theorem \ref{bij1};
    \item flattened Catalan words of length
$n$
and compositions of all even natural numbers into $n-1$ parts of
at most two where the part $0$ is allowed, see Theorem~\ref{thm:bij to comps of even};
\item flattened Catalan words of length $n$ with $k$ runs of weak ascents and binary words  of length $n-1$ where $2k-2$ symbols are replaced with a dot $\bullet$, see Theorem~\ref{thm:bij binary with bullet};
\item flattened Catalan words of length $n$ and Dyck paths of semilength $n$  with $k$ occurrences of $\texttt{DDUU}$, where the height sequence  of occurrences $\texttt{DDU}$ (from left to right) is nondecreasing, see Remark~\ref{rem:flat cat to paths}. 
\item flattened Catalan words of length $n$ and ordered trees with $n$ edges and with $k+1$ nodes having only children as leaves and satisfying two additional conditions, see Remark \ref{rem:flat cat with valleys to trees}.
\end{itemize}

We aggregate our results and the notation used throughout in Table~\ref{tab:notation}.
\begin{table}[H]
\centering
\resizebox{\textwidth}{!}{
\begin{tabular}{|l|c|c|c|c|c|c|}\cline{2-7}
 \multicolumn{1}{l|}{} 	& \multicolumn{6}{c|}{Statistics} \\ \cline{2-7}
 \multicolumn{1}{l|}{} 		& runs  of asc. 		& runs  of w. asc. &runs  of desc. 	&  runs  of w. desc.&$\ell$-valleys	& short valleys   \\ \hline
 Statistic on $w$ & $\runs(w)$ & $\wruns(w)$ &	$\druns(w)$&$\wdruns(w)$ & $\lval(w)$ & 1-$\val(w)$  \\ 
Bivariate g. function & $R(x,y)$ & $W(x,y)$ &$\bar{R}(x,y)$& $\bar{W}(x,y)$&$V_\ell(x,y)$ & $V_1(x,y)$ \\
Distribution & $r(n,k)$ & $w(n,k)$ &$\bar{r}(n,k)$&$\bar{w}(n,k)$& $v_\ell(n,k)$ & $v_1(n,k)$ \\
Total occurrences  over $\f(\Cat_n)$ & $r(n)$  &  $w(n)$ &$\bar{r}(n)$ &$\bar{w}(n)$& $v_\ell(n)$ & $v_1(n)$  \\
\hline
\multicolumn{1}{l|}{} 	&valleys &sym. valleys & $\ell$-peaks & short peaks & peaks & sym.  peaks  \\
\hline
Statistic on $w$&$\val(w)$& $\vsym(w)$& $\lpea(w)$ & 1-$\peak(w)$& $\peak(w)$ & $\psym(w)$  \\
Bivariate g. function &$V(x,y)$& $S(x,y)$ & $P_\ell(x,y)$  & $P_1(x,y)$ & $P(x,y)$ & $T(x,y)$ \\
Distribution &  $v(n,k)$&$s(n,k)$ & $p_\ell(n,k)$ &  $p_1(n,k)$ & $p(n,k)$&$t(n,k)$\\
Total occurrences  over $\f(\Cat_n)$&  $v(n)$& $s(n)$  &  $p_\ell(n)$ & $p_1(n)$ & $p(n)$ & $t(n)$ \\
\hline
\end{tabular}
}
\\
\smallskip
\resizebox{\textwidth}{!}{
\begin{tabular}{|l||c|l|l|}
  \hline
\mbox{Statistic} & \mbox{Bivariate g. f.} & \mbox{Total occurrences over $\f(\Cat_n)$}& \mbox{OEIS}\\
  \hline
  $\runs$  &$\frac{x y (1 - x - x y)}{1 - 2 x + x^2 - 2 x y + x^2 y + x^2 y^2}$&$\frac{1}{4}(3^{n - 1}+1) (n+1) $&\mbox{Not in OEIS}\\[5pt]
   \hline
   $\wruns$ &$\frac{(1 - 2 x)xy}{1 - 4 x + 4 x^2 - x^2 y}$ &$ \frac{1}{36}\left(27 - 9n +(5+n)3^n \right)$&\mbox{Not in OEIS}\\[5pt]
   \hline
   $\druns$ & $\frac{xy(1-2xy)}{1 - 4 x y - x^2 y + 4 x^2 y^2}$ &$\frac{1}{36}\left(27n - 9 +(5n+1)3^n \right)$&\mbox{Not in OEIS}\\[5pt]
   \hline
   $\wdruns $&${\frac {yx \left( 1-xy-x \right) }{{x}^{2}{y}^{2}+{x}^{2}y+{x}^{2}-2
\,xy-2\,x+1}}$
 &$\frac{n+1}{4}(1+3^{n-1}) $& \mbox{Not in OEIS}\\[5pt]
   \hline
   $\lval$&$\frac{x (1 - 2 x + x^{\ell+1} - x^{\ell+1} y)}{(1-x)(1 - 3 x + x^{\ell+1} - x^{\ell+1} y)}$&$\frac{1}{4}\left(1 - 3^{n - 2 - \ell} + 2\cdot 3^{n-2\ell}(n-2-\ell) \right)$&\mbox{Not in OEIS}\\[5pt]
   \hline
   $\val$&$\frac{x - 3 x^2 + x^3 (3 - y)}{(1 - x) (1 - 4 x + 4x^2  - x^2y)}$&$\frac{1}{36}\left(3^n(n-4) + 9n\right)$&\mbox{\seqnum{A212337}}\\[5pt]
   \hline
   $\vsym$&$\frac{x (1 - 2 x) (1 - 2 x + 2 x^2 - x^2 y)}{(1 - x) (1 - 5 x + 8 x^2 - 5 x^3 - x^2 y + 2 x^3 y)}$&$\frac{1}{144}\left(3^n (2 n - 5) - 18 n^2 + 54 n - 27 \right)$&\mbox{Not in OEIS}\\[5pt]
   \hline 
   $\lpea$&$\frac{x(1 - 2 x)}{(1 - x) (1 - 3 x + x^{\ell+1} (1 - y))}$&$\frac{1}{4}\left((3^{n -\ell - 2} (2 n + 1 - 2 \ell)) - 1\right)$&\mbox{Not in OEIS}\\[5pt]
   \hline
   $\peak$&$\frac{x(1-2x)}{1-4x+4x^2-x^2y}$&$\frac14(3^{n-2}-1)(n-1)$&\mbox{\seqnum{A261064}}\\[5pt]
   \hline
   $\psym$&$\frac{x(1-x)(1-2x)}{1 - 5 x + 8 x^2 - 5 x^3 - x^2 y + 2 x^3 y}$&$\frac{1}{144}\left(63 + 3^n + 2 (-45 + 3^n) n + 18 n^2)\right)$&\mbox{Not in OEIS}\\[5pt]
   \hline
\end{tabular}
}
\caption{Summary of notation and results for statistics considered.} \label{tab:notation}
\end{table}

\section{Basic Definitions}\label{sec:Catalan words}

Throughout the article, we will  use the following decomposition of Catalan words, called \emph{first return decomposition} of a Catalan word $w$, which~is \[w=\texttt{0}(w'+1)w'',\] where $w'$ and $w''$ are Catalan words ($w'$ and $w''$ could be empty), and where ($w'+1$) is the word obtained from $w'$ by adding $1$ at all these symbols. Note that whenever $w'$ is the empty word, denoted by $\epsilon$, then $(w'+1)$ remains the empty word.

For example, the first return decomposition of $w=\tx{0122200122322334544}\in \f(\Cat_{19})$ is given by setting $w'=\texttt{0111}$ and $w''=\texttt{00122322334544}$.  
For this word $w$, we have
$\runs(w)=11$, $\wruns(w)=4$, $\druns(w)=16$, $\wdruns(w)=9$, 1-$\val(w)=0$, 2-$\val(w)=2$, $\lval(w)=0$ $(\ell>2)$, $\vsym(w)=1$, 1-$\peak(w)=2$, 2-$\peak(w)=0$, 3-$\peak(w)=1$,   $\lpea(w)=0$  $(\ell>3)$, and $\psym(w)=2$. 

\begin{figure}[H]
\centering
\includegraphics[scale=0.7]{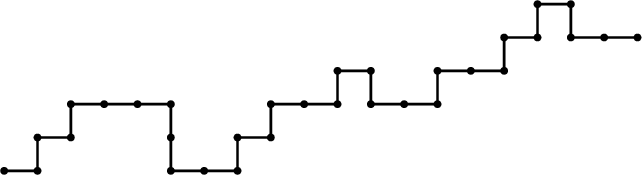}
\caption{Lattice diagram of the word $w=\tx{0122200122322334544}$.} \label{fig1b}
\end{figure}

Drawing Catalan words as lattice diagrams  on the plane proves to be a convenient representation. These diagrams are constructed using unit up steps $(0, 1)$, down steps $(0,-1)$, and horizontal steps $(1,0)$. Each symbol $w_i$ of a Catalan word is represented by the horizontal segment between the points $(i-1, w_i)$ and $(i, w_i)$, and the vertical steps are inserted to obtain a connected diagram.   For example, in Figure~\ref{fig1b}, we illustrate the lattice diagram associated to the Catalan word $w$. 

\begin{remark}
    Let $\Cat^{\uparrow}_n$ denote the set of weakly increasing Catalan words of length $n$. Notice that 
    $|\Cat^{\uparrow}_0|=1$ and  for $n\geq 1$ 
$|\Cat^{\uparrow}_n|=2^{n-1}$, then its generating functions is $1+x/(1-2x)$ if we include the empty word.    Note that the set of nonempty weakly increasing Catalan words is precisely the set of flattened Catalan words with a single weak run. Hence, the generating functions for the later set is $x/(1-2x)$.
\end{remark}

\section{The Distribution of Runs}
\subsection{Runs of Ascents}\label{subsec:runs}
In order to count nonempty flattened Catalan  words according to the length and the number runs of ascents, we introduce the following bivariate generating function
$$R(x,y)=\sum_{w \in\f(\Cat^+)}x^{|w|}y^{\runs(w)}=\sum_{n\geq 1}x^{|w|}\sum_{w\in\f(\Cat_n)}y^{\runs(w)},$$ where the coefficient of $x^ny^k$ is the number of flattened Catalan words of length $n$ with $k$ runs of  ascents.

In Theorem~\ref{thm:fcatruns}, we give an expression for this generating function, but first we provide an example.
\begin{example}
    Consider the flattened Catalan word $w=\texttt{012230123122}\in \f(\Cat_{12})$.
    Then $w$ has $5$ runs of ascents:
        \texttt{012}, \texttt{23}, \texttt{0123},  \texttt{12}, and \texttt{2}.
\end{example}

\begin{theorem}\label{thm:fcatruns}
The generating function for nonempty flattened Catalan words with respect to the length and the number of runs of ascents is
$$R(x,y)=\frac{x y (1 - x - x y)}{1 - 2 x + x^2 - 2 x y + x^2 y + x^2 y^2}.$$
\end{theorem}
\begin{proof}
Let $w$ be a nonempty flattened Catalan word and let $w=\texttt{0}(w'+1)w''$ be the first return decomposition, with $w', w''\in \f(\Cat)$.    There are four different types of this word. Figure \ref{deco2b} illustrates this case. 
          \begin{figure}[H]
 \centering
\includegraphics[scale=0.8]{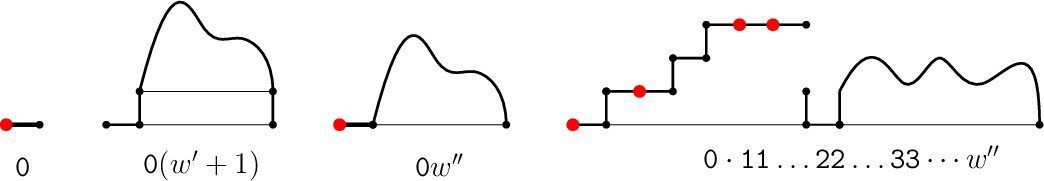}
 \caption{Decomposition of a nonempty flattened Catalan word in $\f(\Cat)$.} \label{deco2b}
 \end{figure}
    
If $w' = w''= \epsilon$, then $w=0$. 
    Then its generating function is $xy$.
    
If $w''=\epsilon$ and $w'\neq \epsilon$, then $w=\texttt{0}(w'+1)$.     Then the generating function is $xR(x,y)$.  

If $w'=\epsilon$ and $w''\neq \epsilon$, then $w=\texttt{0}w''$.
    Then the generating function is $xyR(x,y)$ because we have  an extra run.    
    
If $w'\neq \epsilon$ and $w''\neq \epsilon$, then $w=\texttt{0}(w'+1)w''$.
    Note $w'$ is a weakly increasing word because $w \in \f(\Cat^+)$. 
    Then the bivariate generating function for such words $w'$ is $$\sum_{n\geq1}\sum_{k=1}^{n}\binom{n-1}{k-1}x^ny^k=\sum_{n\geq 0}y(1+y)^{n-1}x^n=\frac{xy}{1-x(1+y)}.$$ 
    Therefore, the generating function for this case is given by $$\frac{x^2y}{1-x-xy}R(x,y).$$ 
   
Therefore, we have the functional equation
$$R(x,y)=xy + x(1+y)R(x,y)+\frac{x^2y}{1-x-xy}R(x,y).$$
Solving this equation, we obtain the desired result. 
\end{proof}

\begin{coro}\label{nber} The generating function for nonempty flattened Catalan words is given by 
$$R(x,1)=\sum_{n\geq 1}\wfwn(n)x^n=\frac{x(1-2x)}{(1-3x)(1-x)}.$$
Therefore, 
$$\wfwn(n)=\frac{1}{2}\left(3^{n-1} + 1\right).$$
\end{coro}

The first few values of the sequence $\wfwn(n)$ ($n\geq 1$) 
correspond to the OEIS entry  \cite[\seqnum{A007051}]{OEIS}:
$$1, \quad 2, \quad 5, \quad 14, \quad 41, \quad 122, \quad 365, \quad 1094, \quad 3281, \quad9842,\ldots.$$

This sequence also counts the compositions of all even natural  numbers (from $0$ to $2(n-1)$) into $n-1$ parts of at most two (the part $0$ is allowed).
\begin{theorem}\label{thm:bij to comps of even}
    Flattened Catalan words of length $n$ and compositions of all even natural numbers (from $0$ to $2(n-1)$) into $n-1$ parts of at most two (the part $0$ is allowed) are in bijection.
\end{theorem}
\begin{proof}
    A bijection $\psi$ between flattened Catalan words of length $n$ and this combinatorial class is given by $\psi(\tx{0})=\epsilon$; $\psi(\tx{0}(w+1))=2\psi(w)$;  $\psi(\tx{0}w)=\tx{0}\psi(w)$; and $\psi(\tx{0}(w+1)w')=1 \psi(w) 1 \psi(w')$.
\end{proof}

\begin{figure}[H]
\centering
 \includegraphics[scale=0.75]{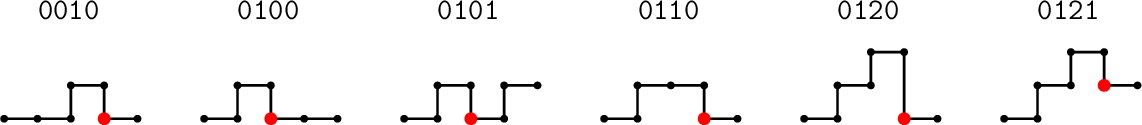}
\caption{Flattened Catalan words of length 4 with 2 runs of ascents. The red marked vertex denotes the start of the second run of ascents.} \label{RunEx1}
\end{figure}

Let $r(n,k)$ denote the number of flattened Catalan words of length $n$ with exactly $k$  runs of  ascents, that is $r(n,k)=[x^ny^k]R(x,y)$, which denotes the coefficient of $x^ny^k$ in $R(x,y)$. The first few rows of this array are 
$$\mathcal{R}:=[r(n,k)]_{n, k\geq 1}=
\begin{pmatrix}
 1 & 0 & 0 & 0 & 0 & 0 & 0 & 0 \\
 1 & 1 & 0 & 0 & 0 & 0 & 0 & 0 \\
 1 & 3 & 1 & 0 & 0 & 0 & 0 & 0 \\
 1 & 6 & \framebox{\textbf{6}} & 1 & 0 & 0 & 0 & 0 \\
 1 & 10 & 19 & 10 & 1 & 0 & 0 & 0 \\
 1 & 15 & 45 & 45 & 15 & 1 & 0 & 0 \\
 1 & 21 & 90 & 141 & 90 & 21 & 1 & 0 \\
 1 & 28 & 161 & 357 & 357 & 161 & 28 & 1
\end{pmatrix}.$$
For example, $r(4,3)=6$, the entry boxed in $\mathcal{R}$ above, and the corresponding flattened Catalan words (and lattice diagrams) are shown in Figure~\ref{RunEx1}.

The array $\mathcal{R}$ 
corresponds to the OEIS entry \cite[\seqnum{A056241}]{OEIS}.
Notice that this sequence has a different combinatorial interpretation. It counts the number of $k$-part order-consecutive partitions of $n$. An {\it order-consecutive partition} of  $\{1,2,\ldots, n\}$ with $k$ parts is a $k$-uplet $(S_1,S_2,\ldots, S_k)$ of subsets such that $S_i\cap S_j=\emptyset$ if $i\neq j$, $\bigcup\limits_{i=1}^k S_i=\{1,2,\ldots, n\}$, where every subset  $S_i$ are in increasing order relatively to their maximum elements, and satisfying the property: for $j=1,\ldots, k$, $\bigcup\limits_{i=1}^{j}S_i$ is an interval (cf. \cite{WM}). 

\begin{theorem}\label{bij1}
Flattened Catalan words of length $n$ with exactly $k$ runs of ascents are in bijection with $k$-part order-consecutive partitions of $n$.
\end{theorem}
\begin{proof}   
We define recursively a map $\psi$ from the set of words in $\f(\Cat_n)$ and the set $\mathcal{OCP}_n$ of order-consecutive partitions of $\{1,2,\ldots , n\}$. 
We consider the four cases of Figure~\ref{deco2b}.
\begin{itemize}
    \item[-] If $w$ belongs to the case ($i$), then $w=\tx{0}$ and we set $\psi(w)=\{1\}$; 
    \item[-] If $w$ belongs to the case ($ii$), then $w=\tx{0}(w'+1)$ and $\psi(w)$ is obtained from $\psi(w')$ by inserting $n$ in the last part; for instance, if $f(w')=\{2,3\}\{1,4\}$, then $f(w)=\{2,3\}\{1,4,5\}$; 
    \item[-] If $w$ belongs to the case ($iii$), then $w=\tx{0}w'$ and $\psi(w)$ is obtained from $\psi(w')$ by adding the part $\{n\}$ on the right; for instance, if $f(w')=\{2,3\}\{1,4\}$, then $f(w)=\{2,3\}\{1,4\}\{5\}$;
    \item[-] 
    If $w$ belongs to the case ($iv$), then $w=w'w''$ where $w'$ consists of one weak run starting with $\tx{01}$. Using the previous cases,  $\psi(w')=S_1\ldots S_k$ where $S_k=\{a_1,\ldots a_\ell,|w'|-1,|w'|\}$ ends with a part containing both $|w'|-1$ and $|w'|$. So, we  set $\psi(w)=S_1\ldots S_{k-1} (\psi(w'')+|w'|-1)\{a_1,\ldots, a_\ell,|w'|-1,|w'|+|w''|\}$. For instance if $w=\tx{0112}~\tx{0120}$, $w'=\tx{0112}$, $w''=\tx{0120}$ and $f(w')=\{1,2\} \{3,4\}$ and  $f(w'')=\{3\}\{1,2,4\}$ then $f(w)=\{1,2\}\{6\}\{4,5,7\}\{3,8\}$. \qedhere
\end{itemize}
\end{proof}
Theorem \ref{bij1}  and \cite[Theorem 6]{WM} imply the following combinatorial expression.
\begin{coro} If $n,k\geq 1$, then 
$$r(n,k)=\sum_{j=0}^{k-1}\binom{n - 1}{2 k - j - 2} \binom{2 k - j - 2}{j}.$$    
\end{coro}

Let $r(n)$ be the total number of runs of ascents over all 
flattened Catalan words  
of length~$n$.
\begin{coro}
We have
$$\sum_{n\geq 0}r(n)x^n=\frac{x - 5 x^2 + 8 x^3 - 3 x^4}{(1 - 3 x)^2 (1 - x)^2}.$$
Moreover, for $n\geq 1$, we have 
$$r(n)=\frac{1}{4}(3^{n - 1}+1) (n+1).$$
\end{coro}

The first few values of the sequence $r(n)$ ($n\geq 1$) are
$$1, \quad  3,\quad 10, \quad 35, \quad 123, \quad 427, \quad 1460, \quad 4923, \quad 16405, \quad 54131,\ldots.$$
This sequence does not appear in the OEIS.

\subsection{Runs of Weak Ascents} In order to count nonempty flattened Catalan  words according to the length and the number runs of weak ascents, we introduce the following bivariate generating function $$W(x,y)=\sum_{w \in \f( \Cat^+)}x^{|w|}y^{\wruns(w)}=\sum_{n\geq 1}x^{|w|}\sum_{w\in\f( \Cat_n)}y^{\wruns(w)},$$ where the coefficient of $x^ny^k$ is the number of flattened Catalan words of length $n$ with $k$ runs of weak ascents. 

\begin{example}
    Consider the flattened Catalan word $w=\tx{012230123122}\in \f(\Cat_{12})$.
Then $w$ has $3$ runs of weak ascents:
        \texttt{01223}, \texttt{0123}, \texttt{122}.
\end{example}

In Theorem~\ref{thm:weakruns}, we give an expression for this generating function.

\begin{theorem}\label{thm:weakruns}
The generating function for the number of nonempty flattened Catalan words with respect to the length and the number of runs of weak ascents is
$$W(x,y)=\frac{(1 - 2 x)xy}{1 - 4 x + 4 x^2 - x^2 y}.$$
\end{theorem}

\begin{proof}
    Let $w$ be a nonempty flattened Catalan word and let $w=\tx{0}(w'+1)w''$ be the
    first return decomposition, with $w', w''\in \f(\Cat)$.    
    There are four different types of this word.
If $w' = w''= \epsilon$, then $w=\tx{0}$. 
    Then its generating function is $xy$.
    If $w''=\epsilon$ and $w'\neq \epsilon$, then $w=\tx{0}(w'+1)$. 
    Then the generating function is $xW(x,y)$.  
    Similarly,  if $w'=\epsilon$ and $w''\neq \epsilon$, then $w=\tx{0}w''$.
    Then the generating function is $xW(x,y)$. 
    If $w'\neq \epsilon$ and $w''\neq \epsilon$, then $w=\texttt{0}(w'+1)w''$.
    Note $w'$ is a weakly increasing
    word because $w \in \f(\Cat^+)$. Then the generating function is given by
    \[x\sum_{k\geq 1}2^kx^kyW(x,y)=\frac{x^2y}{1-2x}W(x,y).\]
         
Therefore, we have the functional equation
$$W(x,y)=xy + 2xW(x,y)+\frac{x^2y}{1-2x}W(x,y).$$
Solving this equation, we obtain the desired result. 
\end{proof}

\begin{figure}[H]
\centering
\includegraphics[scale=0.75]{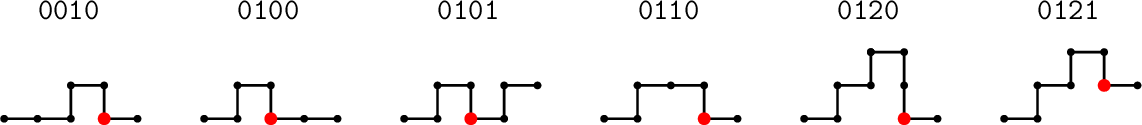}
\caption{Flattened Catalan words of length 4 with 2 runs of weak ascents. The red marked vertex denotes the start of the second run of weak ascents.} \label{RunEx2}
\end{figure}

Let $w(n,k)$ denote the number of flattened Catalan words of length $n$ with exactly $k$ runs of weak ascents, that is $w(n,k)=[x^ny^k]W(x,y)${, which denotes the coefficient of $x^ny^k$ in $W(x,y)$.}
The first few values of this array are 
$$\mathcal{W}:=[w(n,k)]_{n, k\geq 1}=
\begin{pmatrix}
 1 & 0 & 0 & 0 & 0 \\
 2 & 0 & 0 & 0 & 0 \\
 4 & 1 & 0 & 0 & 0 \\
 8 & \framebox{\textbf{6}} & 0 & 0 & 0 \\
 16 & 24 & 1 & 0 & 0 \\
 32 & 80 & 10 & 0 & 0 \\
 64 & 240 & 60 & 1 & 0 \\
 128 & 672 & 280 & 14 & 0 \\
 256 & 1792 & 1120 & 112 & 1
\end{pmatrix}.$$
For example, $w(4,2)=6$, the entry boxed in $\mathcal{W}$ above, and the corresponding flattened Catalan words (and lattice diagrams) are shown in Figure~\ref{RunEx2}. The array $\mathcal{W}$ does not appear in the OEIS.

  \begin{coro}\label{coro:count wnk}
    For $n, k\geq 1$, we have
    $$w(n,k)=2^{n - 2 k + 1}\binom{n-1}{2k-2} .$$
\end{coro}
\begin{proof} From Theorem \ref{thm:weakruns}, we obtain the recurrence relation 
$$w(n,k)-4w(n-1,k)+4w(n-2,k)-4w(n-2,k-1)=0, \quad n\geq 3, k\geq 1,$$
with the initial values $w(2,1)=2$, $w(1,1)=1$, and $w(n,k)$ for $n<k$. It is not difficult to  verify that $2^{n - 2 k + 1}\binom{n-1}{2k-2}$ satisfies the same recurrence relation and the same initial values. Therefore, the sequences are the same. 
\end{proof}
We give an alternate proof of Corollary \ref{coro:count wnk} through a bijective proof. We state the result formally for ease of reference.
\begin{theorem}\label{thm:bij binary with bullet}
    Flattened Catalan words of length $n$ with $k$ runs of weak ascents and binary words  of length $n-1$ where $2k-2$ symbols are replaced with a dot $\bullet$ are in bijection.
\end{theorem}
\begin{proof}
We now give bijection between flattened Catalan words of length $n$ with $k$ runs of weak ascents and binary words  of length $n-1$ where $2k-2$ symbols are replaced with a dot $\bullet$ (Corollary~\ref{coro:count wnk} and a simple combinatorial argument prove that the two classes of objects have the same cardinality).
Let $u=u_1u_2\cdots u_{n-1}$ be such a binary word with $2k-2$ $\bullet$'s, and let us suppose that the $\bullet$'s are on the positions $\{i_1, i_2, \ldots , i_{2k-2}\}$. Then, we define the flattened Catalan words with $k$ runs of weak ascents as follows:

Let $v=v_0v_1\cdots v_{n-1}$ be the word of length $n$ constructed from $u$ by fixing $v_0=\tx{0}$, $v_{i_{2a+1}}:=\tx{1}$, $v_{i_{2a}}:=\tx{0}$, $a=0,1,\ldots, k-1$, and $v_i:=u_i$ for all other positions $i$. We fix $i_0=0$ and $i_{2k-1}=n$.
Now, $v$ consists of the juxtaposition of $k$ nonempty factors of the form $r_a=v_{i_{2a}}\cdots v_{i_{2a+2}-1}$, $a=0,1,\ldots, k-1$, all of them starting with $\tx{0}$. We associate to each factor $s=0s_2\cdots s_p$ the nondecreasing Catalan word $c(s)=\tx{0}c_2\cdots c_{|s|}$, where $c_i=c_{i-1}$ if $s_i=0$ and $c_i=c_{i-1}+\tx{1}$, otherwise (for instance, if $s=011010110$ then $c(s)=\tx{012233455}$).

The bijection $f$ is defined as follows:
$$f(u)=c(r_0)(a_0+c(r_1))(a_0+a_1+c(r_2))\cdots (a_0+a_1+\cdots + a_{k-2}+c(r_{k-1})),$$ where $a_j$ is the number of 
$\tx{1}$'s  in the factor $v_{i_{2(j+1)}}\cdots v_{i_{2(j+1)+1}-1}$.

For instance, if $n=29$ and $k=4$ and $u=10100\bullet 1010\bullet 0110\bullet 01\bullet 0110\bullet 0\bullet 00$. We have $$v=010100\textbf{1}1010~\textbf{0}0110\textbf{1}01~\textbf{0}   0110\textbf{1}0~\textbf{0}00,$$ and \[f(u)=\tx{01122234455}~\tx{22344556}~\tx{4456677}~\tx{666}. \qedhere\]
\end{proof}

Let $w(n)$ be the total number of runs of weak ascents over all flattened Catalan words of length~$n$.

\begin{coro}
For $n\geq 1$, we have
$$\sum_{n\geq 1}w(n)x^n=\frac{x (1 - 2 x)^3}{(1 - 4 x + 3 x^2)^2}.$$
Moreover, for $n\geq 1$, we have 
$$w(n)=\frac{1}{36}\left(27 - 9n +(5+n)3^n \right).$$
\end{coro}

The first few values of the sequence $w(n)$ $(n\geq 1)$ are 
$$1, \quad  2, \quad 6, \quad 20,\quad  67,\quad  222, \quad 728, \quad 2368, \quad 7653, \quad 24602,\ldots$$
This sequence does not appear in the OEIS.

\subsection{Runs of Descents}
In order to count nonempty flattened Catalan  words according to the length and the number runs of descents, we introduce the following bivariate generating function
$$\bar{R}(x,y)=\sum_{w \in \f(\Cat^+)}x^{|w|}y^{\druns(w)}=\sum_{n\geq 1}x^{|w|}\sum_{w\in \f(\Cat_n)}y^{\druns(w)},$$ where the coefficient of $x^ny^k$ is the number of flattened Catalan words of length $n$ with $k$ runs of  descents. 

\begin{example}
    Consider the flattened Catalan word  $w=\tx{012230123122}\in \f(\Cat_{12})$.
    Then $w$ has 10 runs of descents: 
        \texttt{0}, \texttt{1}, \texttt{2}, \texttt{2}, \texttt{30}, \texttt{1}, \texttt{2}, \texttt{31}, \texttt{2},  and \texttt{2}.
\end{example}

It is worth noticing that in any flattened Catalan word $w$ of length $n$, we have $\druns(w)=n+1-\wruns(w)$. Therefore, we can directly deduce Theorem~\ref{teod1} and Corollary~\ref{cordd}. 

\begin{theorem}\label{teod1}
The generating function for the number of nonempty flattened Catalan words with respect to the length and the number of runs of descents is
$$\bar{R}(x,y)=yW\left(xy,\frac{1}{y}\right)=\frac{xy(1-2xy)}{1 - 4 x y - x^2 y + 4 x^2 y^2}.$$
\end{theorem}

\begin{figure}[H]
    \centering
\includegraphics[scale=0.75]{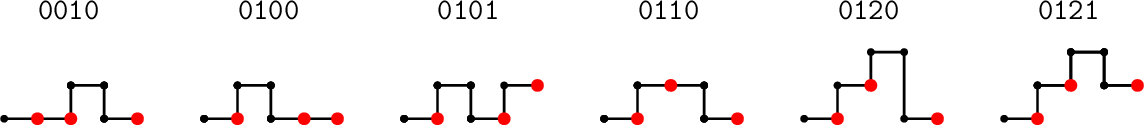}
    \caption{Flattened Catalan words of length 4 with 3 runs of descents.  The red marked vertices denote the end  of a  run of descents. }\label{RunDescentsEx}
\end{figure}

 Let $\bar{r}(n,k)$ denote the number of flattened Catalan words of length $n$ with exactly $k$  runs of  descents, that is $\bar{r}(n,k)=[x^ny^k]\bar{R}(x,y)$, which denotes the coefficient of $x^ny^k$ in $\bar{R}(x,y)$.
 The first few values of this arrays are 
 $$\bar{\mathcal{R}}\coloneq[\bar{r}(n,k)]_{n, k\geq 1}=
 \begin{pmatrix}
 1 & 0 & 0 & 0 & 0 & 0 & 0 & 0 & 0 \\
 0 & 2 & 0 & 0 & 0 & 0 & 0 & 0 & 0 \\
 0 & 1 & 4 & 0 & 0 & 0 & 0 & 0 & 0 \\
 0 & 0 & \framebox{\textbf{6}} & 8 & 0 & 0 & 0 & 0 & 0 \\
 0 & 0 & 1 & 24 & 16 & 0 & 0 & 0 & 0 \\
 0 & 0 & 0 & 10 & 80 & 32 & 0 & 0 & 0 \\
 0 & 0 & 0 & 1 & 60 & 240 & 64 & 0 & 0 \\
 0 & 0 & 0 & 0 & 14 & 280 & 672 & 128 & 0 \\
 \end{pmatrix}.$$
 For example, $\bar{r}(4,3)=6$, the entry boxed in $\bar{\mathcal{R}}$ above, and the corresponding flattened Catalan words (and lattice diagrams) are shown in Figure~\ref{RunDescentsEx}.
 The array $\bar{\mathcal{R}}$ does not appear in the OEIS.

\begin{coro}\label{cordd}
     For $n, k\geq 1$, we have
     $$\bar{r}(n,k)=2^{2k-n-1}\binom{n - 1}{2 (n - k) }.$$
\end{coro}
A combinatorial interpretation of this last formula can be obtained from the bijection $f$ (see Section~3.2) between flattened Catalan words of length $n$ with $n+1-k$ runs of weak ascents (or equivalently with $k$ descents) and binary words of length $n-1$ with  $(2n-2k)$ dots $\bullet$.

Let $\bar{r}(n)$ be the total number of runs of descents over all 
flattened Catalan words  
of length~$n$.
\begin{coro}
We have
$$\sum_{n\geq 0}\bar{r}(n)x^n=\frac{x (1 - 4 x + 4 x^2 + 2 x^3)}{(1 - 4 x + 3 x^2)^2}.$$
Moreover, for $n\geq 1$, we have
$$\bar{r}(n)=\frac{1}{36}\left(27n - 9 +(5n+1)3^n \right).$$
\end{coro}

The first few values of the sequence $\bar{r}(n)$ ($n\geq 1$) are
$$1, \quad 4, \quad 14, \quad 50,\quad  179,\quad  632, \quad 2192, \quad 7478, \quad 25157, \quad 83660,\ldots.$$
This sequence does not appear in the OEIS.

\subsection{Runs of Weak Descents}
In a flattened Catalan word of length $n$, the number of runs of ascents plus the number of runs of weak descents equals $n+1$. Hence, the number $\bar{w}(n,k)$ of flattened Catalan words of length $n$ with $k$ runs of weak descents equals the number $r(n,k)$ of flattened Catalan words of length $n$  with $k$ runs of ascents. Moreover, we can defined a simple involution $\phi$ on  $\f(\Cat_n)$ such that $\phi(w)=w'$ with $\wdruns(\phi(w))=\runs(w)$, as follows: $\phi(\epsilon)=\epsilon$, $\phi(\tx{0}(w+1))=\tx{0}\phi(w)$,   $\phi(\tx{0}w)=\tx{0}(1+\phi(w))$, and $\phi(\tx{0}(1+w)w')=\tx{0}(1+\phi(w))\phi(w')$ whenever $w,w'\neq \epsilon$.
Then, we the results can be restated as those in Section~\ref{subsec:runs}. 
\begin{theorem}
The generating function for the number of nonempty flattened Catalan words with respect to the length and the number of runs of weak descents is
$$\bar{W}(x,y)=R(x,y)={\frac {yx \left(1- xy-x \right) }{{x}^{2}{y}^{2}+{x}^{2}y+{x}^{2}-2
\,xy-2\,x+1}}.$$
Therefore, $$\bar{w}(n,k)=r(n,k)=\sum_{j=0}^{k-1}\binom{n - 1}{2 k - j - 2} \binom{2 k - j - 2}{j}.$$ 
\end{theorem}

\begin{coro}
We have
$$\sum_{n\geq 0}\bar{w}(n)x^n=\sum_{n\geq 0}r(n)x^n={\frac {x \left( 1-3\,{x}^{3}+8\,{x}^{2}-5\,x \right) }{ \left( 3\,{x
}^{2}-4\,x+1 \right) ^{2}}}.$$
Moreover, for $n\geq 1$, we have $$\bar{w}(n)=r(n)=\frac{n+1}{4}(1+3^{n-1}).$$ 
\end{coro}

\section{The Distribution of Valleys}

\subsection{Valleys}

In order to count nonempty flattened Catalan words according to the length and the number $\ell$-valleys, we introduce the following bivariate generating function $$V_\ell(x,y)=\sum_{w \in\f(\Cat^+)}x^{|w|}y^{\lval(w)}=\sum_{n\geq 1}x^{|w|}\sum_{w\in\f(\Cat_n)}y^{\lval(w)},$$
where $\lval(w)$ denotes the number of occurrences of subwords of the form $ab^\ell(b+1)$, and $a> b$, in $w$. 
The coefficient of $x^ny^k$ in $V_\ell (x,y)$ is the number of flattened Catalan words of length $n$ with $k$ $\ell$-valleys.

In Theorem~\ref{teosv1l}, we give an expression for this generating function.

\begin{theorem}\label{teosv1l}
The generating function for nonempty flattened Catalan words with respect to the length and the  number of $\ell$-valleys is
$$V_\ell(x,y)=\frac{x (1 - 2 x + x^{\ell+1} - x^{\ell+1} y)}{(1-x)(1 - 3 x + x^{\ell+1} - x^{\ell+1} y)}.$$
\end{theorem}

\begin{proof} Let $w$ be a nonempty flattened Catalan word, and let  $w=\texttt{0}(w'+1)w''$ be the first return decomposition, with $w', w''\in \f(\Cat)$. 
If $w'=w''=\epsilon$, then $w=\texttt{0}$, and its generating function is $x$.
If $w'\neq \epsilon$ and $w''=\epsilon$, then $w=\texttt{0}(w'+1)$, and its generating function is $xV_\ell (x,y)$.
Similarly,  if $w'=\epsilon$ and $w''\neq \epsilon$, then $w=\texttt{0}w''$, and its generating function is $xV_\ell (x,y)$.
Finally, if $w'\neq \epsilon$ and $w'' \neq \epsilon$, then $w=\texttt{0}(w'+1)w''$.
Because $w$ is a flattened Catalan word, $w'$ must be a 
weakly increasing 
word,
and we distinguish two cases. 
If $w''$ is of the form $\tx{0}^{\ell-1}w'''$, where $w'''$ starts with $\tx{01}$, 
then $w=\texttt{0}(w'+1)\texttt{0}^{\ell-1}w'''$, and the generating function is 
$$\left(\frac{x^{\ell+1}y}{1-2x}\right)\left(V_\ell(x,y)-(x+xV_{\ell}(x,y)\right).$$
Notice that $T_\ell(x,y):=V_\ell(x,y)-(x+xV_{\ell}(x,y))$ is obtained using the complement of the generating function for the word $\texttt{0}$ and the words starting with $\texttt{00}$.

The second case is the negation, so, $w''$ does not start with $\texttt{0}^\ell\texttt{1}$. Notice that $\ell $ is fixed because we are interested in the $\ell$-valleys, so  the generating function is 
$$\frac{x^2}{1-2x}(V_\ell(x,y)-x^{\ell-1} T_\ell(x,y)).$$
Therefore, we have the functional equation
\begin{align*}
    V_{\ell}(x,y) &= x+2xV_{\ell}(x,y)+
\left(\frac{x^{\ell+1}y}{1-2x}\right)T_\ell(x,y) +
 \frac{x^2}{1-2x}(V_\ell(x,y)-x^{\ell-1} T_\ell(x,y)).
\end{align*}
Solving this equation, we obtain the desired result.
\end{proof}

\begin{figure}[H]
\centering
\includegraphics[scale=0.75]{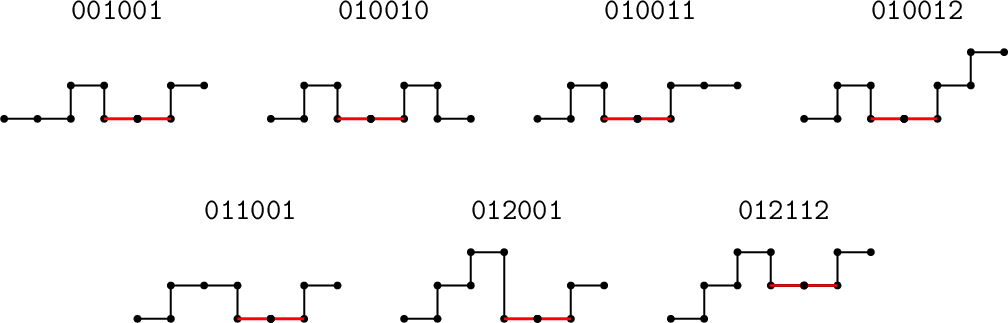}
\caption{Flattened Catalan words of length 6 with one $2$-valley. The red edges indicate the location of the $2$-valley.} \label{flatCatExlVal}
\end{figure}

Let $v_\ell(n,k)$ denote the number of flattened Catalan words of length $n$ with exactly $k$ $\ell$-valleys, that is $v_\ell(n,k)=[x^ny^k]V_\ell(x,y)$, which denotes the coefficient of $x^ny^k$ in $V_\ell(x,y)$.
For example, the first few values of this array for $\ell=2$ are
$$\mathcal{V}_2\coloneq[v_2(n,k)]_{n\geq 4, k\geq 0}=
\begin{pmatrix}
 14&0&0&0\\
 40&1&0&0\\
 115&\framebox{\bf{7}}&0&0\\
 331&34&0&0\\
 953&140&1&0\\
 2744&527&10&0\\
 7901&1877&64&0
    \end{pmatrix}.$$
For example, $v_2(6,1)=7$, the entry boxed in $\mathcal{V}_2$ above, and the corresponding flattened Catalan words of length $6$ with one $2$-valley (and lattice diagrams) are shown in Figure~\ref{flatCatExlVal}.
The first column of the array $\mathcal{V}_2$ corresponds to OEIS entry 
\cite[\seqnum{A052963}]{OEIS}.

 Let $v_\ell(n)$ be the sum of all $\ell$-valleys in the set of flattened Catalan words of length $n$. 
 \begin{coro}
 The generating function of the sequence
 $v_\ell(n)$  is 
 $$\sum_{n\geq 1}v_\ell(n)x^n=\frac{x^{\ell+3}}{(1-x)(1-3x)^2}.$$
 Moreover, for $n\geq 1$, we have
 \begin{align*}
     v_\ell(n)=\frac{1}{4}\left(1 - 3^{n - 2 - \ell} + 2\cdot 3^{n-2\ell}(n-2-\ell) \right).
 \end{align*}
 \end{coro}

Taking $\ell=1$ in Theorem~\ref{teosv1l}, we obtain the generating function for nonempty flattened Catalan words with respect to the length and the  number of short valleys
 $$V_1(x,y)=\sum_{w \in\f(\Cat^+)}x^{|w|}y^{\text{1-}\val(w)}=\frac{x - 2 x^2 + x^3 (1 - y)}{(1 - x) (1 - 3 x + x^2 (1 - y))}.
 $$

\begin{figure}[H]
\centering
\includegraphics[scale=0.8]{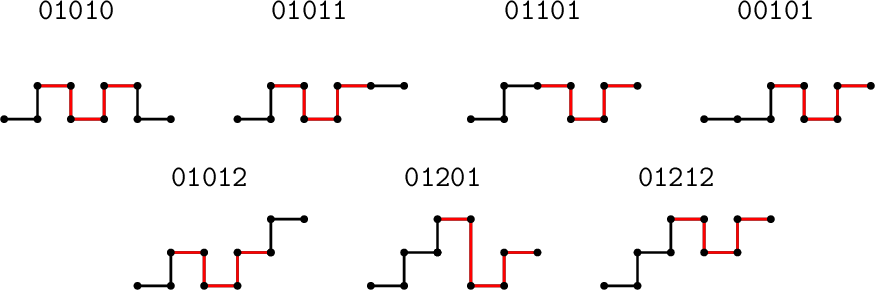}
\caption{Flattened Catalan words of length 5 with one  short valley. The red edges indicates the location of the  short valley.} \label{n5oneshortvalley}
\end{figure}

 Let $v_1(n,k)$ denote the number of flattened Catalan words of length $n$ with exactly $k$ short valleys, that is $v_1(n,k)=[x^ny^k]V_1(x,y)$, which denotes the coefficient of $x^ny^k$ in $V_1(x,y)$. The first few values of this array are
 $$\mathcal{V}_1=[v_1(n,k)]_{n\geq 1, k\geq 0}=
 \begin{pmatrix}
  1 & 0 & 0 & 0 \\
 2 & 0 & 0 & 0 \\
 5 & 0 & 0 & 0 \\
 13 & 1 & 0 & 0 \\
 34 & \framebox{\textbf{7}} & 0 & 0 \\
 89 & 32 & 1 & 0 \\
 233 & 122 & 10 & 0 \\
 610 & 422 & 61 & 1 \\
 1597 & 1376 & 295 & 13 
 \end{pmatrix}.$$
 For example, $v_1(5,1)=7$, the entry boxed in $\mathcal{V}_1$ above, and the corresponding flattened Catalan words of length $5$ with exactly one short valley (and lattice diagrams) are shown in Figure~\ref{n5oneshortvalley}.

\begin{remark}\label{rem:flat cat to paths}
In \cite{Paper1}, we proved that Catalan words of length $n$ with $k$ short valleys are in one-to-one correspondence with Dyck paths of semilength $n$ with $k$ occurrences of $\texttt{DDUU}$. Taking the restriction  on flattened Catalan words of this bijection, we obtain a one-to-one correspondence between flattened Catalan words of length $n$ and Dyck paths of semilength $n$  with $k$ occurrences of $\texttt{DDUU}$, where the height sequence  of occurrences $\texttt{DDU}$ (from left to right) is nondecreasing. 
\end{remark}

We can also obtain the generating function for the
 number of flattened Catalan words of length $n$ with respect to the number of valleys (we consider all $\ell$-valleys for $\ell\geq 1$). 

 \begin{theorem}\label{teov1}
 The generating function for nonempty flattened Catalan words with respect to the length and the  number of valleys is
 $$V(x,y)=\frac{x - 3 x^2 + x^3 (3 - y)}{(1 - x) (1 - 4 x + 4x^2  - x^2y)}.$$
 \end{theorem}

 Let $v(n,k)$ denote the number of flattened Catalan words of length $n$ with exactly $k$ valleys, that is $v(n,k)=[x^ny^k]V(x,y)$, which denotes the coefficient of $x^ny^k$ in $V(x,y)$.
 The first few values of this arrays are 
 $$\mathcal{V}=[v(n,k)]_{n\geq 1, k\geq 0}=
 \begin{pmatrix}
  1 & 0 & 0 & 0 \\
 2 & 0 & 0 & 0 \\
 5 & 0 & 0 & 0 \\
 13 & 1 & 0 & 0 \\
 33 & 8 & 0 & 0 \\
 81 & 40 & 1 & 0 \\
 193 & 160 & \framebox{\textbf{12}} & 0 \\
 449 & 560 & 84 & 1 \\
 1025 & 1792 & 448 & 16 
 \end{pmatrix}.$$
 For example, $v(7,2)=12$, the entry boxed in $\mathcal{V}$ above, and the corresponding flattened Catalan words of length $7$ with exactly two valleys are \begin{align*}
 \texttt{0010101}, \quad 
  \texttt{0100101}, \quad
   \texttt{0101001}, \quad
 \texttt{0101010}, \quad
 \texttt{0101011}, \quad
 \texttt{0101012},\\
  \texttt{0101101}, \quad
 \texttt{0101201}, \quad
  \texttt{0101212}, \quad   
 \texttt{0110101},\quad
 \texttt{0120101},\quad
 \texttt{0121212}.
 \end{align*}

 \begin{coro}
     For $n\geq 0$ we have 
      $$v(n,k)=\begin{cases} (n - 1)2^{n - 2} + 1, & \text{if } k=0\\
          2^{n-2k-2}\binom{n-1}{2 k +1}, & \text{if } k\geq 1
      \end{cases}.$$
\end{coro}
Note that $v(n,0)$ corresponds to OEIS entry \cite[\seqnum{A005183}]{OEIS}.

\begin{remark}\label{rem:flat cat with valleys to trees}
In \cite{Paper1}, we proved that Catalan words of length $n$ with $k$  valleys are in one-to-one correspondence with ordered trees with $n$ edges and having exactly $k+1$ nodes all of those children are leaves. Taking the restriction  on flattened Catalan words of this bijection, we obtain a one-to-one correspondence between flattened Catalan words of length $n$ and ordered trees with $n$ edges and with $k+1$ nodes having only children as leaves and satisfying the following:
\begin{itemize}
    \item if $T_1, T_2, \ldots, T_r$ are the subtrees of the root, then $T_i$, $i\in[1,r-1]$, is  nondecreasing (i.e. for any node, its subtrees, except the rightmost,  consist of one node only), 
\item the rightmost subtree of the root again satisfies all these properties. 
\end{itemize}
\end{remark}

 Let $v(n)$ be the sum of all valleys in the set of flattened  Catalan words of length $n$. 
 \begin{coro}
 The generating function of the sequence
 $v(n)$  is 
 $$\sum_{n\geq 0}v(n)x^n=\frac{x^4}{(1-x)^2(1-3x)^2}.$$ 
 Moreover, for $n\geq 4$, we have 
\begin{align*}
     v(n)=\frac{1}{36}\left(3^n(n-4) + 9n\right).
 \end{align*}
\end{coro}

For $n\geq 4$, the first few values of the sequence $v(n)$ are 
$$
	1,\quad  8,\quad  42,\quad  184,\quad  731,\quad  2736,\quad  9844,\quad  34448,\quad  118101,\quad  398584,\ldots.$$
This sequence corresponds to OEIS entry \cite[\seqnum{A212337}]{OEIS}.

\subsection{Symmetric Valleys}
A \emph{symmetric  valley} is a valley of the form $a(a-1)^\ell a$ with $\ell\geq 1$. Let $\vsym(w)$ denote  the number of  symmetric valleys in the word $w$. 
In order to count flattened Catalan words according to the length and the number of symmetric valleys, we introduce the following bivariate generating function
generating function
$$S(x,y)=\sum_{w\in\f(\Cat^+)}x^{|w|}y^{\vsym(w)}=\sum_{n\geq 1}x^{|w|}\sum_{w\in\f(\Cat_n)}y^{\vsym(w)},$$
where the coefficient of $x^ny^k$ in $S(x,y)$ is the number of nonempty flattened Catalan words of length~$n$ with $k$ symmetric $\ell$-valleys.

In Theorem~\ref{symval}, we give an expression for this generating function.

\begin{theorem}\label{symval}
    The generating function of the nonempty flattened Catalan words with respect to the length and the number of symmetric valleys is 
    \[S(x,y)=\frac{x (1 - 2 x) (1 - 2 x + 2 x^2 - x^2 y)}{(1 - x) (1 - 5 x + 8 x^2 - 5 x^3 - x^2 y + 2 x^3 y)}.\]  
\end{theorem}

\begin{proof} Let $w$ be a nonempty flattened Catalan word, and let  $w=\texttt{0}(w'+1)w''$ be the first return decomposition, with $w', w''\in \Cat$. 
If $w'=w''=\epsilon$, then $w=\texttt{0}$, and its generating function is $x$. 
If $w'\neq \epsilon$ and $w''=\epsilon$, then $w=\texttt{0}(w'+1)$, and its generating function is $xS(x,y)$. 
Similarly, if $w'=\epsilon$ and $w''\neq \epsilon$, then $w=\texttt{0}w''$, and its generating function is $xS(x,y)$.  
Finally, if $w'\neq \epsilon$ and $w''\neq \epsilon$, then $w=\texttt{0}(w'+1)w''$, we consider three cases.
\begin{enumerate}
    \item
    If $w'=\tx{0}^k$ and $w''$ has a nonzero entry,
     then its generating function is
    \[\left(\frac{x^2}{1-x}\right)y\left(S(x,y)-\frac{x}{1-x}\right).\]
    \item If $w'$ is a weakly increasing flattened Catalan word different than $\tx{0}^k$, and $w''$ has a nonzero entry,
     then its generating function is
    \[x\left( \frac{x}{1-2x}- \frac{x}{1-x}\right)\left(S(x,y)-\frac{x}{1-x}\right).\]
    \item If $w'$ is a weakly increasing flattened Catalan word and $w''=\tx{0}^{k}$,
    then its generating function~is
    \[\frac{x^3}{(1-x)(1-2x)}.\]
    \end{enumerate}
  Therefore, we have the functional equation
\begin{multline*}
    S(x,y)=x+2xS(x,y)+\left(\frac{x^2}{1-x}\right)y\left(S(x,y)-\frac{x}{1-x}\right)+\\
x\left( \frac{x}{1-2x}- \frac{x}{1-x}\right)\left(S(x,y)-\frac{x}{1-x}\right) + \frac{x^3}{(1-x)(1-2x)}.
\end{multline*}
Solving the obtained functional equation yields the desired result.
\end{proof}

Let $s(n,k)$ denote the number of flattened Catalan words of length $n$ with exactly $k$ symmetric valleys, that is $s(n,k)=[x^ny^k]S(x,y)$, which denotes the coefficient of $x^ny^k$ in $S(x,y)$.
The first few values of this arrays are 
$$\mathcal{S}=[s(n,k)]_{n\geq 1, k\geq 0}=
\begin{pmatrix}
 1 & 0 & 0 & 0 & 0 \\
 2 & 0 & 0 & 0 & 0 \\
 5 & 0 & 0 & 0 & 0 \\
 13 & 1 & 0 & 0 & 0 \\
 34 & \framebox{\textbf{7}} & 0 & 0 & 0 \\
 90 & 31 & 1 & 0 & 0 \\
 242 & 113 & 10 & 0 & 0 \\
 659 & 375 & 59 & 1 & 0 
\end{pmatrix}.$$

\begin{figure}[H]
\centering
\includegraphics[scale=0.75]{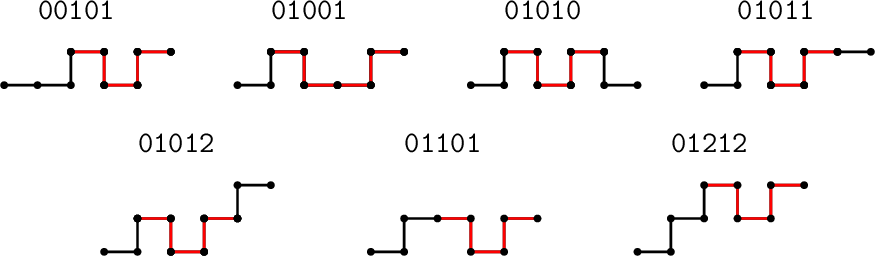}
\caption{Flattened Catalan words of length 5 with one symmetric valley. In red we mark the location of the symmetric valley.} \label{flatCatExSymVal}
\end{figure}

For example, $s(5,1)=7$, the entry boxed in $\mathcal{S}$ above, and the corresponding flattened Catalan words of length 5 with 1 symmetric valley are given in Figure \ref{flatCatExSymVal}. The array $\mathcal{S}$ does not appear in the OEIS.

Let $s(n)$ be the sum of all symmetric valleys in the set of flattened Catalan words of length $n$. 
 \begin{coro}\label{corosym}
 The generating function of the sequence
 $s(n)$  is 
$$\sum_{n\geq 0}s(n)x^n=\frac{x^4 (1 + 2 x)}{(1 - 3 x)^2 (1 - x)^3}.$$
 Moreover, for $n\geq 4$, we have
 \begin{align*}
     s(n)=\frac{1}{144}\left(3^n (2 n - 5) - 18 n^2 + 54 n - 27 \right).
 \end{align*}
 \end{coro}

The first few values of the sequence $s(n)$ $(n\geq 4)$ are
$$ 
 1, \quad 7, \quad 33, \quad 133, \quad 496,
 \quad 1770, \quad 6142, \quad 20902,
 \quad 70107,\quad 232489, \dots.$$
This sequence does not appear in the OEIS.

\section{The Distribution of Peaks}

\subsection{Peaks}
In order to count flattened Catalan words according to the length and the number of $\ell$-peaks, we introduce the following bivariate generating function
$$P_\ell(x,y)=\sum_{w \in\f(\Cat^+)}x^{|w|}y^{\lpea(w)}=\sum_{n\geq 1}x^{|w|}\sum_{w\in\f(\Cat_n)}y^{\lpea(w)},$$
where $\lpea(w)$ denotes the number of occurrences of subwords of the form $a(a+1)^\ell b$, and $a\geq b$, in $w$.  
The coefficient of $x^ny^k$ in $P_\ell (x,y)$ is the number of flattened Catalan words of length $n$ with $k$ $\ell$-peaks.

In Theorem~\ref{teosp1l}, we give an expression for this generating function.
\begin{theorem}\label{teosp1l}
The generating function for nonempty flattened Catalan words with respect to the length and the  number of $\ell$-peaks is
\[P_{\ell}(x,y)=\frac{x(1 - 2 x)}{(1 - x) (1 - 3 x + x^{\ell+1} (1 - y))}.
\]
\end{theorem}
\begin{proof} 
Let $w$ be a nonempty flattened Catalan word, and let  $w=\texttt{0}(w'+1)w''$ be the first return decomposition, with $w', w''\in \Cat$. 
If $w'=w''=\epsilon$, then $w=\texttt{0}$, and its generating function is $x$. 
If $w'\neq \epsilon$ and $w''=\epsilon$, then $w=\texttt{0}(w'+1)$, and its generating function is $xP_{\ell}(x,y)$. 
Similarly, if $w'=\epsilon$ and $w''\neq \epsilon$, then $w=\texttt{0}w''$, and its generating function is $xP_{\ell}(x,y)$.  
Finally, if $w'\neq \epsilon$ and $w''\neq \epsilon$, then $w=\texttt{0}(w'+1)w''$, its generating function is 
$$x\left(\frac{x}{1-2x}-x^\ell- \frac{x^{\ell+1}}{1-2x}\right)P_{\ell}(x,y)+xy\left(x^\ell+\frac{x^{\ell+1}}{1-2x}\right)P_{\ell}(x,y).$$
Therefore, we have the functional equation
\begin{align*}
P_{\ell}(x,y)&=x+2xP_{\ell}(x,y)+x\left(\frac{x}{1-2x}-x^\ell-\frac{x^{\ell+1}}{1-2x}\right)P_{\ell}(x,y)\\&\hspace{.5in}+xy\left(x^\ell+\frac{x^{\ell+1}}{1-2x}\right)P_{\ell}(x,y).
\end{align*}

Solving the obtained functional equation yields the desired results.
\end{proof}

Let $p_\ell(n)$ be the sum of all $\ell$-peaks in the set of flattened Catalan words of length $n$. 
\begin{coro}
The generating function of the sequence
$p_\ell(n)$  is 
$$\sum_{n\geq 1}p_\ell(n)x^n=\frac{x^{\ell+2}(1 - 2 x)}{(1 - 3 x)^2 (1 - x)}.$$
Moreover, for $n\geq 1$ we have
\begin{align*}
p_\ell(n)=\frac{1}{4}\left((3^{n -\ell - 2} (2 n + 1 - 2 \ell)) - 1\right).
\end{align*}
\end{coro}

Taking $\ell=1$ in Theorem~\ref{teosp1l}, establishes that the generating function for flattened Catalan words with respect to the length and the number of short peaks is
$$P_1(x,y)=\frac{x(1-2x)}{(1 - x) (1 - 3 x + x^2 (1 - y))}.
$$
Let $p_1(n,k)$ denote the number of flattened Catalan words of length $n$ with exactly $k$ short peaks, that is $p_1(n,k)=[x^ny^k]P_1(x,y)$, which denotes the coefficient of $x^ny^k$ in $P_1(x,y)$. The first few values of this array are
$$\mathcal{P}_1=[p_1(n,k)]_{n\geq 1, k\geq 0}=
\begin{pmatrix}
 1 & 0 & 0 & 0 &0\\
 2 & 0 & 0 & 0 &0\\
 4 & 1 & 0 & 0& 0\\
 9 & 5 & 0 & 0 &0\\
 22 & 18 & 1 & 0&0 \\
  56 & 58 & \framebox{\textbf{8}} & 0 &0\\
 145 & 178 & 41 & 1 &0\\
 378 & 532 & 173 & 11 &0\\
 988 & 1563 & 656 & 73&1 
\end{pmatrix}.$$
For example, $p_1(6,2)=8$, the entry boxed in $\mathcal{S}$ above, and the corresponding flattened Catalan words of length 6 with 2 short peaks are
\begin{align*}
\texttt{001010},\,
\texttt{010100},\,
\texttt{010101},\,
\texttt{010010},\,
\texttt{010120},\,
\texttt{010121},\,
\texttt{012010},\,
\texttt{012121}.
\end{align*}

While the full array $\mathcal{P}_1$ does not appear in the OEIS, for $n\geq 1$ we have 
$p_1(n,0)=F_{2(n-1)}+1$, where $F_m$ is the $m$th Fibonacci number with initial values $F_1=F_2=1$. For $n\geq1$, the sequence $p_1(n,0)$ corresponds to the OEIS entry \cite[\seqnum{A055588}]{OEIS}.

Using a similar proof as for Theorem~\ref{teosp1l}, we generalize the result in order to obtain the following generating function for the
 number of flattened Catalan words of length $n$ with respect to the number of peaks (we consider all $\ell$-peaks for $\ell\geq 1$).

\begin{theorem}\label{teop1}
The generating function for flattened Catalan words with respect to the length and the  number of peaks is
$$P(x,y)=\frac{x(1-2x)}{1-4x+4x^2-x^2y}.$$
\end{theorem}

Let $p(n,k)$ denote the number of flattened Catalan words of length $n$ with exactly $k$ peaks, that is $p(n,k)=[x^ny^k]P(x,y)$, which denotes the coefficient of $x^ny^k$ in $P(x,y)$.
The first few values of this arrays are 
 $$\mathcal{P}=[p(n,k)]_{n\geq 1, k\geq 0}=
 \begin{pmatrix} 
  1 & 0 & 0 & 0 & 0 \\
  2 & 0 & 0 & 0 & 0 \\
  4 & 1 & 0 & 0 & 0 \\
  8 &  \framebox{\textbf{6}} & 0 & 0 & 0 \\
  16 & 24 & 1 & 0 & 0 \\
  32 & 80 & 10 & 0 & 0 \\
  64 & 240 & 60 & 1 & 0 \\
 128 & 672 & 280 & 14 & 0 \\
 256 & 1792 & 1120 & 112 & 1 \\
 \end{pmatrix}.
$$
For example, $p(4,1)=6$, the entry boxed in $\mathcal{P}$ above, and the corresponding flattened Catalan words of length 4 with 1 peaks are
\begin{align*}
\texttt{0010},\quad
\texttt{0100},\quad
\texttt{0110},\quad
\texttt{0101},\quad
\texttt{0120},\quad
\texttt{0121}.
\end{align*}
The array $\mathcal{P}$  does not appear in the OEIS. 

 Let $p(n)$ be the sum of all peaks in the set of flattened Catalan words of length $n$. 
\begin{coro}
The generating function of the sequence
$p(n)$  is 
$$\sum_{n\geq 0}p(n)x^n=\frac{(1 - 2 x) x^3}{(1 - 4 x + 3 x^2)^2}.$$
Moreover, for $n\geq 3$, we have 
\[p(n)=\frac14(3^{n-2}-1)(n-1).\]
\end{coro}

The first few values of the sequence $p(n)$ ($n\geq 3$) are 
 $$1, \quad 6, \quad 26, \quad 100, \quad 363, \quad 1274, \quad  4372, \quad  14760,\quad  14760, \quad 49205,\dots .$$
This sequence corresponds to the OEIS entry \cite[\seqnum{A261064}]{OEIS}. Our combinatorial interpretation is new.
%\pamela{This is the "Number of non-selfintersecting broken lines in a convex $(n+1)$-gon", should we try ot give a bijection or no?}\jluc{I don't see how!:-(}

\subsection{Symmetric Peaks}
A \emph{symmetric peak} is a peak of the form $a(a+1)^\ell a$ with $\ell\geq 1$.  Let $\psym(w)$ denote  the number of the symmetric peaks of the word $w$.  In order to count flattened Catalan words according to the length and the number symmetric peaks, we introduce the following bivariate generating function
$$T(x,y)=\sum_{w\in\f(\Cat^+)}x^{|w|}y^{\psym(w)}=\sum_{n\geq 1}x^{|w|}\sum_{w\in\f(\Cat_n)}y^{\psym(w)},$$
where the coefficient of $x^ny^k$ in $T(x, y)$ is the number of flattened Catalan words of length~$n$ with $k$ symmetric peaks.

Theorem~\ref{thm:sympeaks}, we give an expression for this generating function.
 \begin{theorem}\label{thm:sympeaks}
     The generating function of the nonempty flattened Catalan words with respect to the length and the number of symmetric peaks is 
\[T(x,y)=\frac{x(1-x)(1-2x)}{1 - 5 x + 8 x^2 - 5 x^3 - x^2 y + 2 x^3 y}.\]
\end{theorem}
\begin{proof} 
Let $w$ be a nonempty flattened Catalan word, and let  $w=\texttt{0}(w'+1)w''$ be the first return decomposition, with $w', w''\in \f(\Cat)$.
If $w'=w''=\epsilon$, then $w=\texttt{0}$, and its generating function is $x$. 
If $w'\neq \epsilon$ and $w''=\epsilon$, then $w=\texttt{0}(w'+1)$, and its generating function is $xT(x,y)$. 
Similarly,  if $w'=\epsilon$ and $w''\neq \epsilon$, then $w=\texttt{0}w''$, and its generating function is $xT(x,y)$.

Finally, if $w'\neq \epsilon$ and $w''\neq \epsilon$, then $w=\texttt{0}(w'+1)w''$, and we have two cases to consider.
\begin{enumerate}
    \item If $w'$ is all $\tx{0}$'s, its generating function is
    \[\frac{x^2y}{1-x}T(x,y).\]
    \item 
    Otherwise, the generating function is
    \[x\left(\frac{x}{1-2x}-\frac{x}{1-x}\right)T(x,y).\]
\end{enumerate}
Therefore, we have the functional equation is
\[T(x,y)=x+2xT(x,y)+\frac{x^2y}{1-x}T(x,y)+x\left(\frac{x}{1-2x}-\frac{x}{1-x}\right)T(x,y).\]
Solving this equation yields the desired result.
\end{proof}

Let $t(n,k)$ denote the number of flattened Catalan words of length $n$ with exactly $k$ symmetric peaks, that is $t(n,k)=[x^ny^k]T(x,y)$, which denotes the coefficient of $x^ny^k$ in $T(x,y)$.
 The first few values of this arrays are 
 $$\mathcal{T}=[t(n,k)]_{n\geq 1, k\geq 0}=\begin{pmatrix}
 1 & 0 & 0 & 0 & 0 \\
 2 & 0 & 0 & 0 & 0 \\
 4 & 1 & 0 & 0 & 0 \\
 9 & \framebox{\textbf{5}}  & 0 & 0 & 0 \\
 23 & 17 & 1 & 0 & 0 \\
 63 & 51 & 8 & 0 & 0 \\
 176 & 149 & 39 & 1 & 0 \\
 491 & 439 & 153 & 11 & 0  
\end{pmatrix}.
$$

\begin{figure}[H]
    \centering
    \includegraphics[scale=0.75]{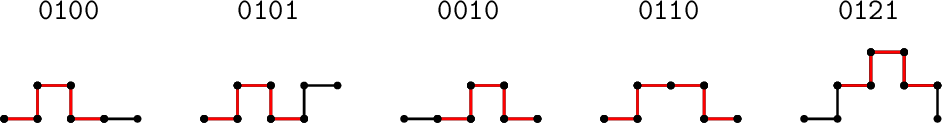}
\caption{Flattened Catalan words of length 4 with 1 symmetric peak. In red we mark the location of the symmetric peak.}\label{SymPeaksEx}
\end{figure}

 For example, $t(4,1)=5$, the entry boxed in $\mathcal{T}$ above, and the corresponding flattened Catalan words of length 4 with 1 symmetric peak (and lattice diagrams) are shown in Figure~\ref{SymPeaksEx}.
The first and second column of the array ${\mathcal{T}}$ coincides with OEIS entries 
\cite[\seqnum{A369328}, \seqnum{A290900}]{OEIS}.
  The full array ${\mathcal{T}}$ does not appear in the OEIS.

Let $t(n)$ be the sum of all symmetric peaks in the set of flattened Catalan words of length $n$. 
\begin{coro}
 The generating function of the sequence
 $t(n)$  is 
 $$\sum_{n\geq 0}t(n)x^n=\frac{(1 - 2 x)^2 x^3}{(1 - 3 x)^2 (1 - x)^3}.$$
  Moreover, for $n\geq 3$,  we have
   \begin{align*}
 t(n)&=\frac{1}{144}\left(63 + 3^n + 2 (-45 + 3^n) n + 18 n^2)\right).
 \end{align*}
 \end{coro}

For $n\geq 3$, the first few values of the sequence $t(n)$ are 
\[1,\quad 5,\quad 19,\quad 67,\quad 230,\quad 778,\quad 2602,\quad 8618,\quad 28303, \quad 92275,\ldots .\]
This sequence does not appear in the OEIS.

\bigskip
 \noindent{\bf Acknowledgement:}
Jean-Luc Baril was  supported by University of Burgundy. Pamela E.~Harris was supported in part by a Karen Uhlenbeck EDGE Fellowship.
Jos\'e L.~Ramírez was partially supported by Universidad Nacional de Colombia.  The authors thank Kimberly J. Harry and Matt McClinton  for their helpful discussions during the completion of this manuscript.

\end{document}